\documentclass{article}

\usepackage{epsfig,psfrag}
\usepackage{amsmath}
\usepackage{amssymb}
\usepackage{amsthm}
\usepackage{hyperref}

\newtheorem{theorem}{Theorem}

\begin{document}
\title{ Throwing a Ball as Far as Possible, Revisited}
\author{Joshua Cooper \\ Anton Swifton \\ Department of Mathematics \\ University of South Carolina}
\markright{Throwing a Ball}
\date{\today}

\maketitle

\begin{abstract}
What initial trajectory angle maximizes the arc length of an ideal projectile? We show the optimal angle, which depends neither on the initial speed nor on the acceleration of gravity, is the solution \(\theta\) to a surprising transcendental equation: \(\csc(\theta) = \coth(\csc(\theta))\), i.e., \(\theta = \csc^{-1}(\alpha)\) where \(\alpha\) is the unique positive fixed point of \(\coth(x)\).  Numerically, \(\theta \approx 0.9855 \approx 56.47^\circ\).  The derivation involves a nice application of differentiation under the integral sign.
\end{abstract}

Maximizing how far from the origin an ideal projectile lands is a classic problem commonly used as an exercise in introductory physics and calculus courses. The problem of maximizing the length of the projectile's trajectory, however, is not usually considered. We show that the latter problem can be used as a neat application of differentiation under the integral sign, a useful but lesser-known version of Leibniz's rule. Furthermore, the solution satisfies a surprisingly simple transcendental formula.  For present purposes, an ``ideal projectile'' is the path of a point-mass traveling over a flat surface, experiencing constant downward acceleration (\`{a} la gravity), with no resistance or propulsion.

An alternative approach that avoids integration under the integral sign and provides a numerical solution only is described in \cite{sarafian} and \cite{ju}. Those works, however, do not yield the transcendental equation for \(\theta\). Strangely, this equation -- and, of course, the same root \(\theta\) -- occurs implicitly in a different physics problem (see \cite{szirtes}): the angle at which a catenary hangs that minimizes tension at the point of attachment!  We would like to understand why this equation occurs in a seemingly unrelated setting.

\begin{figure}[h] \label{fig:fig1}
\centering
\includegraphics[width=347pt, height=164.5pt]{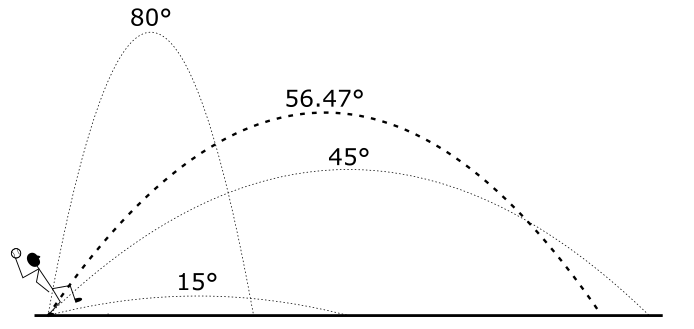}
\caption{Illustration of trajectories for various initial angles.  Note the longest horizontal distance (\(45^\circ\)) and longest arc length (\(\approx 56.47^\circ\)).}
\end{figure}

\begin{theorem} The initial trajectory angle \(\theta\) that maximizes the ground-to-ground arc length traveled by an ideal projectile is the smallest positive root of \(\csc(\theta) = \coth(\csc(\theta))\).  It is also given by \(\theta = \csc^{-1}(\alpha)\), where \(\alpha\) is the unique positive fixed point of \(\coth(x)\).   Numerically, \(\theta \approx 0.9855 \approx 56.47^\circ\).
\end{theorem}

Assume that the projectile is launched with initial speed \(v\) at angle \(\theta\) to the horizon, and that the acceleration of gravity is \(g\).  The length of its trajectory can be found by integrating its speed over time.
\[
L(\theta) = \int\limits_0^{\tau(\theta)} \sqrt{v^2 \cos^2\theta + (v \sin\theta - gt)^2} \, dt.
\]
Here
\[\tau(\theta) = \frac{2v\sin(\theta)}{g}\]
is the time when the projectile hits the ground, i.e., the solution 
\(t = \tau(\theta)\) to \(v t \sin \theta - g t^2/2 = 0\). It will be useful to note before we attempt to maximize \(L(\theta)\) that
\[
v \sin\theta - g \tau(\theta) = -v \sin\theta,
\]
and therefore, 
\[
\sqrt{v^2 \cos^2\theta + (v \sin\theta - g \tau(\theta))^2} = v,
\]
which makes sense due to the symmetry of the picture. The horizontal component of the velocity is constant, and the vertical component at the time when the projectile hits the ground is the negative of its value at the beginning. Therefore, the speed at the end (the left-hand side) is equal to the speed at the beginning (the right-hand side).

Note that \(\theta = \pi/2\) does not yield the longest trajectory, despite being the angle that maximizes {\em vertical distance traveled}. We can calculate
\[
L\left(\frac{\pi}{2}\right) = \int\limits_0^{\tau(\pi/2)}\sqrt{(v-gt)^2} \, dt = \int\limits_0^{2v/g} |v-gt| dt = \frac{v^2}{g}
\]
and estimate 
\[
L\left(\frac{\pi}{4}\right) = \int\limits_0^{\tau(\pi/4)}\sqrt{\frac{v^2}{2} + \left(\frac{v\sqrt{2}}{2} - gt\right)^2}dt > \int\limits_0^{v\sqrt{2}/g} \frac{v}{\sqrt{2}} \, dt = \frac{v^2}{g}.
\]
The inequality is strict due to the fact that the square of the vertical component of the velocity \(\left(\frac{v\sqrt{2}}{2} - gt\right)^2\) is strictly positive on the whole interval except for the apex.

From now on we consider \(\theta \neq \pi/2\). To find the maximum of \(L(\theta)\) we first have to find its critical points. We can calculate the derivative using the technique of differentiation under the integral sign. In \cite{Dieu} the rule of ``differentiation under the integral sign'' is stated as follows (in a more general form).

\begin{theorem}
If \(f(x, \alpha)\), \(a = a(\alpha)\), and \(b = b(\alpha)\) are continuously differentiable with respect to \(\alpha\), then 
\[
\frac{d}{d \alpha} \int\limits_a^b f(x, \alpha) dx = \int\limits_a^b \frac{\partial f}{\partial \alpha} dx + f(b, \alpha) \frac{db}{d\alpha} - f(a, \alpha) \frac{da}{d\alpha}.
\]
\end{theorem}

In our case, \[f(t, \theta) = \sqrt{v^2 \cos^2\theta + (v \sin\theta - gt)^2},\] and \[\frac{\partial f}{\partial \theta} = \frac{1}{2 \sqrt{v^2 \cos^2\theta + (v \sin\theta - gt)^2}} \frac{\partial}{\partial \theta} \left(v^2 \cos^2\theta + (v \sin\theta - gt)^2\right).\] Therefore, \(\partial f/\partial \theta\) is continuous whenever the denominator is positive, which is always true, as long as \(\theta \neq \pi / 2\). Applying differentiation under the integral sign to \(L(\theta)\) yields
\begin{align*}
L'(\theta) &= \frac{d}{d\theta}\int\limits_0^{\tau(\theta)} \sqrt{v^2 \cos^2\theta + (v \sin\theta - gt)^2} \, dt \\
& = \sqrt{v^2 \cos^2\theta + (v \sin\theta - g\tau(\theta))^2} \tau'(\theta) + \int\limits_0^{\tau(\theta)} \frac{\partial}{\partial\theta} \left(\sqrt{v^2 \cos^2\theta + (v \sin\theta - gt)^2}\right) dt\\
& = v \tau'(\theta) + \int\limits_0^{\tau(\theta)} \frac{-2v^2 \cos\theta \sin\theta + 2(v \sin\theta - gt)v \cos\theta}{2 \sqrt{v^2 \cos^2\theta + (v \sin\theta - gt)^2}} \, dt\\
& = v \tau'(\theta) +  v \cos\theta \int\limits_0^{\tau(\theta)} \frac{-gt}{\sqrt{v^2 \cos^2\theta + (v \sin\theta - gt)^2}} \, dt.
\end{align*}
The integral can be computed using simple trigonometric substitution.  Letting \(h(t) = v \sin\theta - gt\) for readability's sake,
\begin{align*}
\int\limits_0^{\tau(\theta)} &\frac{-gt}{\sqrt{v^2 \cos^2\theta + h(t)^2}} \,dt 
\\
&= \int\limits_0^{\tau(\theta)} \frac{h(t)}{\sqrt{v^2 \cos^2\theta + h(t)^2}} \,dt 
- \int\limits_0^{\tau(\theta)} \frac{v \sin\theta}{\sqrt{v^2 \cos^2\theta + h(t)^2}} \,dt
\\
&= -\frac{1}{2g} \int\limits_0^{\tau(\theta)} \frac{d(v^2 \cos^2\theta + h(t)^2)}{\sqrt{v^2 \cos^2\theta + h(t)^2}} + \frac{v\sin\theta}{g} \int\limits_0^{\tau(\theta)} \frac{dh}{\sqrt{v^2 \cos^2\theta + h(t)^2}}
\\
&= \left. -\frac{1}{g} \sqrt{v^2 \cos^2\theta + h(t)^2} \right|_{0}^{\tau(\theta)} 
+ \left. \frac{v\sin\theta}{g} \ln\left|h(t) + \sqrt{h(t)^2 + v^2\cos^2\theta}\right| \right|_{0}^{\tau(\theta)} 
\\
&= \frac{v\sin\theta}{g} \ln \left(\frac{1-\sin\theta}{1 + \sin\theta}\right).
\end{align*}

\begin{figure}
\includegraphics[width=221.75pt, height=165.75pt]{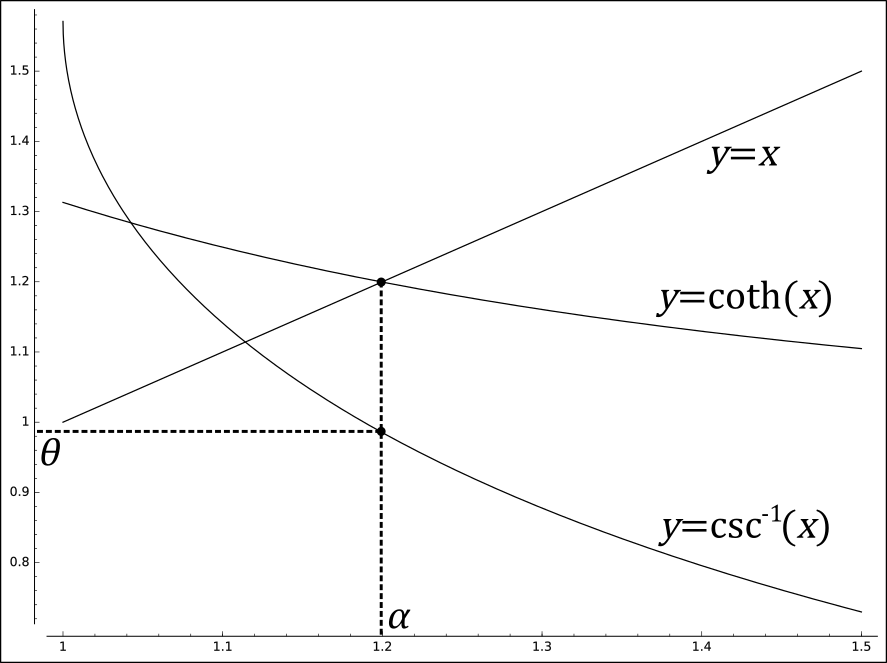}\centering
\caption{Visualizing the fixed point \(\alpha\) of \(\coth\) and its image under \(\csc^{-1}\).}
\end{figure}

Therefore, setting \(L'(\theta)\) equal to zero yields
\[
L'(\theta) = \frac{2v^2\cos\theta}{g} + \frac{v^2\sin\theta\cos\theta}{g} \ln \left(\frac{1-\sin\theta}{1 + \sin\theta}\right) = 0.
\]
Note that \(v^2 > 0\) and \(g > 0\) factor out of this equation, as does \(\cos \theta\) because \(\theta \neq \pi/2\). The resulting equation is
\[
\frac{1}{2} \sin\theta \ln \left(\frac{1+\sin\theta}{1 - \sin\theta}\right) = \sin\theta \tanh^{-1}(\sin \theta) = 1,
\]
which can be further simplified, using the fact that \(\tanh^{-1}(x) = \coth^{-1}(1/x)\), to 
\[\csc\theta = \coth(\csc \theta).\]

Now, \(\coth(x)\) has a unique positive fixed point since it is a continuous decreasing function on \((0, \infty)\), it tends to infinity when \(x\) tends to zero, and to one when \(x\) tends to infinity. So, there if \(x\) is sufficiently small, then \(\coth(x) > x\), and if \(x\) is sufficiently large, \(\coth(x) < x\). Both functions \(\coth(x)\) and \(x\) are continuous and differentiable on \((0,\infty)\), \(x\) is increasing, and \(\coth(x)\) is decreasing; therefore they intersect in exactly one point as a corollary of Rolle's theorem. Thus, if \(\alpha\) is the unique positive fixed point of \(\coth\), then \(\theta = \csc^{-1}( \alpha\)). This equation can be solved numerically (e.g., by the SageMath system \cite{sage}), and has exactly one solution on the interval \([0, \frac{\pi}{2}]\), which is approximately \(0.9855 \approx 56.47^\circ\).  Assuming natural units, i.e., \(g = v = 1\), the length of the longest possible trajectory is \(1.19967864\ldots\), while \(\pi/4\) yields a trajectory of length \(1.14779357\ldots\), a more than \(4.5\%\) difference.

\begin{figure}
\includegraphics[width=175.75pt, height=86.5pt]{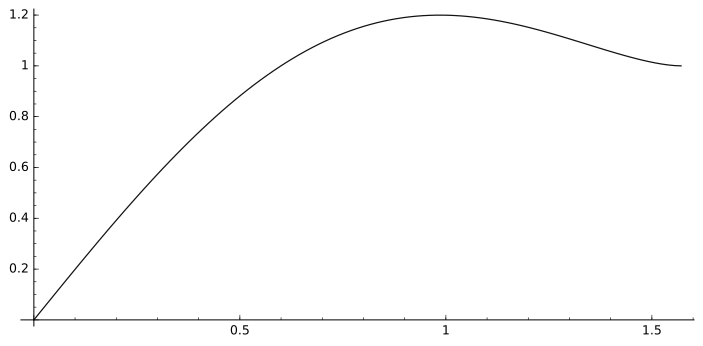}\centering
\caption{Arc length of the trajectory as a function of the angle \(\theta\).}
\end{figure}

\end{document}